\documentclass[journal,twoside,web]{ieeecolor}
\usepackage{lcsys}
\usepackage{cite}
\usepackage{amsmath,amssymb,amsfonts}
\usepackage{algorithmic}
\usepackage{graphicx}   % For including images
\usepackage{caption}    % For customizing captions
\usepackage{subcaption} % For subfigures

\usepackage{amsmath}
\usepackage{amsthm, amssymb}
\theoremstyle{definition}
\newtheorem{definition}{Definition}
\theoremstyle{assumption}

\theoremstyle{problem}

\theoremstyle{lemma}
\newtheorem{lemma}{Lemma}
\theoremstyle{remark}
\newtheorem{remark}{Remark}
\theoremstyle{theorem}
\newtheorem{theorem}{Theorem}

\DeclareMathOperator*{\argmin}{arg\,min}

\usepackage{breqn}

\usepackage{textcomp}
\def\BibTeX{{\rm B\kern-.05em{\sc i\kern-.025em b}\kern-.08em
    T\kern-.1667em\lower.7ex\hbox{E}\kern-.125emX}}
\begin{document}
\title{Trajectory-Oriented Control Using Gradient Descent: An Unconventional Approach}
\author{Ramin Esmzad, Hamidreza Modares, \IEEEmembership{Senior, IEEE}
\thanks{This work is supported by the National Science Foundation (NSF) award ECCS-2227311. }
\thanks{R. Esmzad and H. Modares are with Michigan State University, East Lansing, MI 48824 USA (e-mails:  (esmzadra, modaresh)@msu.edu).}
% \thanks{S. ....}
% \thanks{T. ...}
}

\maketitle

\begin{abstract}
In this work, we introduce a novel gradient descent-based approach for optimizing control systems, leveraging a new representation of stable closed-loop dynamics as a function of two matrices i.e.
 the step size or direction matrix and
 value matrix of the Lyapunov cost function. This formulation provides a new framework for analyzing and designing feedback control laws. We show that any stable closed-loop system can be expressed in this form with appropriate values for the step size and value matrices. Furthermore, we show that this parameterization of the closed-loop system is equivalent to a linear quadratic regulator for appropriately chosen weighting matrices. We also show that trajectories can be shaped using this approach to achieve a desired closed-loop behavior. 
\end{abstract}

% \begin{abstract}
% In this work, we introduce a novel gradient descent-based approach for optimizing control systems, leveraging a new representation of stable closed-loop dynamics as
% $I-2\Gamma P$. Here, $I$ is the identity matrix, 
% $\Gamma$ is a matrix that represents the step size or direction matrix of the gradient descent algorithm, and
% $P$ is a value matrix of the Lyapunov cost function $V_k=x_k^\top P x_k$. This formulation provides a new framework for analyzing and designing feedback control laws. We show that any stable closed-loop system can be expressed in this form with appropriate values of $\Gamma$ and $P$. Furthermore, we show that this parameterization of the closed-loop system with $\Gamma$ is equivalent to a linear quadratic regulator for appropriately chosen $Q$ and $R$ matrices. We also show that trajectories can be shaped using this approach to achieve a desired closed-loop behavior. 
% \end{abstract}

\begin{IEEEkeywords}
Gradient descent, Optimal control, Optimization, Linear Matrix Inequality (LMI).
\end{IEEEkeywords}

\section{Introduction}
\label{sec:introduction}
\IEEEPARstart{G}{radient} descent (GD)  and its variants~\cite{ruder2016overview}, including accelerated optimization algorithms \cite{laborde2020lyapunov}, natural gradient \cite{rattray1998natural,martens2020new}, contraction-based optimization \cite{TSUKAMOTO2021135,wensing2020beyond}, and projected or distributed GD \cite{8186925} have been widely used in the system and control domain. With the recent advancement of machine learning in control theory to achieve autonomy despite uncertainties, GD algorithms play a key paradigm in modern control systems. Often, the process of learning a system model (in model-based control design) or a control policy (in model-free control design) is encapsulated in the form of a parameterized model or controller. These unknown parameters are learned to optimize a cost function that captures the system modeling error or the cost of the control implementation. In this line of work, the GD is used by the controller to achieve some specifications. We call it \textit{GD-based control}.

In model-based and model-free GD-based control methods, the former includes parameter identification for parametric uncertain systems \cite{ljung1998system} or learning the weights of neural networks \cite{80202,sutton2018reinforcement} that are used to predict the behavior of a system with structured uncertainties. The latter includes GD-based learning of an optimal control policy using various reinforcement learning (RL) algorithms, such as value learning \cite{8169685} and policy gradient \cite{7752775} methods, GD-based adaptive control to adjust the control parameters based on the real-time behavior of the system \cite{gaudio2019connections,ioannou2012robust,10599619} or to match the output of a reference model \cite{landau2011adaptive}, GD-based model predictive control (MPC) to optimize the control input by minimizing a cost function over a prediction horizon \cite{bharadhwaj2020model,ASWANI20131216},  and many others (e.g., direct data-driven control \cite{soudbakhsh2023data,HOU20133}, multi-agent control \cite{4118472}, fault detection and isolation \cite{isermann2006fault,9451628}). 
% In fault detection and isolation (FDI), gradient descent is used to minimize residuals or fault indicators, enabling the detection and isolation of faults in control systems. This ensures the reliability and safety of critical systems, such as in aerospace and industrial automation~\cite{isermann2006fault}.

% In data-driven control approaches, gradient descent is utilized to optimize the control policy based on the data directly without explicitly identifying the system model. This is particularly useful in complex systems where modeling is challenging~\cite{HOU20133}.

% Gradient descent is employed in the design of controllers for distributed systems, such as multi-agent systems or networked control systems. The algorithm optimizes local control laws to achieve global performance objectives~\cite{4118472}.

% Gradient descent is utilized in machine learning-based tuning of control systems, where algorithms learn optimal controller settings from data. This method is especially useful in complex systems where manual tuning is difficult, such as in chemical processes or energy systems.~\cite{jesawada2022modelbasedreinforcementlearningapproach}. Gradient descent can be employed to tune PID controllers by minimizing a performance criterion, such as the integral of the squared error (ISE). This approach ensures optimal controller parameters for the desired system response~\cite{LAWRENCE2022105046}.\\

% Gradient descent can be used to optimize controller parameters to achieve robustness against uncertainties and disturbances~\cite{HAUSWIRTH2024100941}.

On the other hand, in another line of work, control tools have been used to improve the performance of GD algorithms \cite{lessard2016analysis,GUNJAL2024100417,padmanabhan2024analysis}. We call it \textit{controlled GD}. Controlled GD methods view GD algorithms through the lens of control theory by treating the optimization process as a dynamic system that can be controlled to improve convergence and stability. To this end, Lyapunov stability \cite{laborde2020lyapunov}, robust control theory \cite{padmanabhan2024analysis}, and passivity immersion \cite{nayyer2022passivity} are well-known control tools that have been leveraged. In sharp contrast to GD-based control, which uses GD for control purposes, the focus of controlled GD is to use control theoretic tools to modify the gradient vector field and improve its convergence. 

In GD-based control approaches, the GD algorithm is used to adjust the parameters of an algorithm to achieve control specifications. The core idea is to change the parameters of the underlying mechanism based on a predefined cost function. A fundamental assumption in all these methods is that an appropriate cost function is available for which its optimization leads to a desired behavior in the state space domain\cite{Cothren2021DataenabledGF}. However, it is a daunting challenge to find or learn an appropriate cost function that captures the designer's intention. Typically, a controller is learned, and its trajectory is observed on a simulated (oracle) or real system, and an appropriate cost is found by trial and error until it leads to a desired behavior. This is also related to the reward hacking problem in the context of RL \cite{skalse2022defining}. Therefore, a fundamental flaw in the learning process of the optimization-based control design using GD is that the behaviors of the system trajectories in the search space cannot be predicted for the optimal controller before its learning and deployment. Besides, a control approach that directly guides the system trajectories can resolve conflicts that arise between multiple objectives or cost functions. 
% For example, even if the cost function of an RL agent is appropriately chosen for a set of system constraints, if the constraints change, the RL trajectories can be in huge conflict with the safety constraints. Frequent safety interventions can cause the RL trajectory to significantly violate its original one. While a new cost function for which the RL agent level sets are more aligned with the safe set can significantly improve the overall performance, there is no mechanism in the standard control approaches to do so. 

As a paradigm shift, this paper presents a different perspective on the concept of control design using the GD algorithm. We present a GD-like closed-loop system dynamics from which the controller can be designed to directly guide trajectories while optimizing a cost function. To this end, we look at the states of the system as parameters of the cost function and use the GD algorithm to shape the trajectories of the system. In other words, we update the trajectories of the system using the GD algorithm and extract a policy based on it. We take advantage of a new representation of stable closed-loop dynamics as 
a deviation of the multiplication of two matrices (step size and value) from the identity matrix. We show that any stable closed-loop system can be expressed in this form, and this parameterization is equivalent to a linear quadratic regulator (LQR) with suitably chosen $Q$ and
$R$ matrices. Simulation results show that the presented approach can shape the trajectories to achieve desired behaviors.

The paper is organized as follows. Section \ref{sec:2} reviews the GD method.
Section \ref{sec:3} formulates the trajectory-oriented control using the GD concept.
Section \ref{sec:4} explores its connection to the LQR.
Section \ref{sec:5} extends the approach to the heavy-ball gradient method.
Section \ref{sec:6} presents the simulation results.
Section \ref{sec:7} concludes the paper.

\section{A review of GD algorithm} \label{sec:2}
This section outlines the standard GD concept for solving optimization problems. We will leverage GD and its insight to establish a trajectory-oriented control approach based on GD-like dynamics. 

\indent GD is a fundamental optimization algorithm widely used to minimize cost functions by iteratively moving toward the steepest descent of the cost function. This method is especially prevalent in machine learning, neural network training, and various numerical optimization problems because of its simplicity and effectiveness. Consider an optimization problem in which the objective is to find the minimum of a differentiable function \( f: \mathbb{R}^n \to \mathbb{R} \). Formally, the problem is stated as
\[
x^* = \argmin_{x \in \mathbb{R}^n} f(x)
\]
where \( x \) is an \( n \) dimensional vector of decision variables and $x^*$ is the optimal solution.
The GD algorithm aims to adjust the decision vector \( x \) iteratively to reduce the value of the objective function \( f \). The core idea is to move in the direction opposite to the gradient of \( f \), which indicates the direction of the steepest ascent. By moving against this direction, the algorithm progresses towards the minimum. The update rule for the GD algorithm at iteration \( k \) is given by
\[
x_{k+1} = x_k - \alpha \nabla f(x_k)
\]
where \( x_k \) is the current point in the parameter space, \( \alpha > 0 \) is the step size (or learning rate), controlling the magnitude of the update, and
 \( \nabla f(x_k) \) is the gradient of the objective function \( f \) evaluated at \( x_k \). The gradient \( \nabla f(x_k) \) is a vector of partial derivatives, given by
\[
\nabla f(x_k) = \left[ \frac{\partial f}{\partial x_{1,k}}, \frac{\partial f}{\partial x_{2,k}}, \ldots, \frac{\partial f}{\partial x_{n,k}} \right]^T,
\] where $x_{i,k}$ is the $i-$th element of $x_k$.
The gradient vector $\nabla f(x_k)$ points in the direction of the steepest increase of \( f \). Hence, the negative gradient \( -\nabla f(x_k) \) points toward the steepest decrease. 
For the GD algorithm to converge to a (possibly local) minimum, the step size \( \alpha \) must be chosen appropriately. The convergence is influenced by several factors such as
\begin{itemize}
    \item Step Size (\( \alpha \)): If \( \alpha \) is too large, the algorithm may overshoot the minimum, causing divergence. If \( \alpha \) is too small, the convergence will be slow. Adaptive methods or line search techniques are often used to dynamically adjust \( \alpha \).
    \item Function Properties: To ensure convergence to a global optimum, the objective function \( f \) should be convex and differentiable. For non-convex functions, the GD may converge to a local minimum.
\end{itemize}

% \subsection{A new interpretation of the Gradient descent algorithm for dynamical systems}
% In this subsection just briefly touch on how we can use the same concept (GD) for shaping states. 
% \\
% What is the gradient dominated cost function?

\section{A new interpretation of the GD algorithm for dynamical systems} \label{sec:3}
% \textcolor{red}{First check this concept for a nonlinear system and see what happens. If its interpretable, add it to the paper, Then describe it for linear systems. }
This section formulates the trajectory-oriented control using GD-like dynamics and designs a static state feedback controller for dynamical systems. It shows how one can directly guide the trajectories of a dynamical system to achieve the desired behavior using the GD concept. 
\subsection{Problem Statement} \vspace{-4pt}
Consider the following deterministic linear time-invariant (LTI) discrete-time system
\begin{align}
    x_{k+1} = A x_k + B u_k, \label{eq:sys}
\end{align}
where $k \in \mathbb{N}$, $x_k \in \mathbb{R}^n$ is the system's state, and $u_k \in \mathbb{R}^m$ is the control input. Moreover, $A$ and $B$ are transition and input matrices of appropriate dimensions, respectively. 

The idea is to design a control policy $u_k$
% \begin{align}
%     u_k = K x_k,
% \end{align}
to force the system states to descend the gradient of a cost function $V_k$. In this paper, we show that this can be achieved by the following trajectory-oriented GD dynamics
\begin{align}
    x_{k+1} = x_k - \Gamma \frac{\partial V_k}{\partial x_k}, \label{eq:gdl}
\end{align}
where  $\Gamma \succ 0$ is the step size matrix that can be used to shape the trajectories of the closed-loop system, and $V_k$ is chosen as 
\begin{align}
    V_k = x^T_k P x_k, \label{eq:Vk}
\end{align}
 for some $P \succ 0$. 
 % The step size matrix $\Gamma$ in general is 
 % \begin{align}
 %     \Gamma = \text{diag}(\gamma_1, \gamma_2, \cdots, \gamma_n), \quad \gamma_i > 0, \,\, \forall \, i=1, \cdots, n.
 % \end{align}
The trajectory-oriented GD dynamics \eqref{eq:gdl} will span the trajectories of a closed-loop system from which a controller will be found. If $\Gamma = \gamma I$, where $\gamma > 0$, on the existence of a $P\succ 0$ that realizes a closed-loop system, the state trajectories will descend in the opposite direction of the gradient of $V_k$. This can be quite helpful since the behavior of the system's trajectories for a closed-loop system can be predicted. Otherwise, for a general $\Gamma$, the trajectories will descend with a different angle (not the exact opposite) that can be tuned to guide trajectories to achieve a desired behavior.
\begin{remark}
    Unlike the prevalent GD methods that are used to update the parameters of the controller or neural networks, here we use the GD concept to update and guide the states of a dynamical system directly. A control law is then found based on this trajectory-oriented GP dynamics with desired behaviors. 
\end{remark}
\begin{definition}\label{Definition.1} \textbf{$\lambda-$Contractive and Positive Invariant Sets}
Consider system \eqref{eq:sys}. Let $\lambda \in (0,1]$. A convex and compact set $\mathcal{E}$ which includes the origin as its interior is considered $\lambda$-contractive if for any $x_k \in \mathcal{E}$, the next state falls into the shrunken version of it i.e. $x_{k+1} \in \lambda \mathcal{E}$. The set $\mathcal{E}$ is named positive invariant whenever $\lambda = 1$. $\lambda$-contractivity implies exponential stability and specifies an exponential decay rate convergence of trajectories.
\end{definition}
% $\lambda$-contractivity is a stronger condition than most forms of stability because it implies exponential convergence of trajectories. While a $\lambda$-contractive system is always stable (specifically, exponentially stable), a stable system is not necessarily $\lambda$-contractive. The key distinction lies in the rate and manner of convergence. $\lambda$-contractivity specifies an exponential decay rate, whereas general stability conditions allow for a wider range of convergence behaviors.
\subsection{Stability Analysis of Trajectory-oriented GD Dynamics}
While the presented trajectory-oriented GD can be leveraged to achieve various desired behaviors, here, we establish the stability of the trajectory-oriented GD system~\eqref{eq:gdl}. The following lemma outlines the conditions for its stability. 
\begin{lemma} \label{lemma:gd}
    The system \eqref{eq:gdl} with the cost function \eqref{eq:Vk} is $\lambda-$contractive with $0< \lambda \leq 1$ (and thus stable) if there exists a $Y = P^{-1}\succ 0 $ and $\Gamma \succ 0$ such that the following linear matrix inequality (LMI) holds 
    \begin{align}
    \begin{bmatrix}
        \lambda Y & (Y-2\Gamma )^\top \\ Y-2\Gamma  & Y
    \end{bmatrix} \succ 0. \label{eq:gdlmi}
\end{align}
\textbf{Proof}. The gradient of the cost $V_k$ with respect to the states is $\frac{\partial V_k}{\partial x_k}=2Px_k$, so the GD system \eqref{eq:gdl} can be simplified as 
\begin{align}
    x_{k+1} = ( I - 2\Gamma P) x_k, \label{eq:gdl2}
\end{align}
which parameterizes the closed-loop system as a product of two matrices $2\Gamma P$ deviated from the identity matrix (we will show this in the next section). 
The rest of the proof is straightforward and is based on Lyapunov's stability analysis. Please refer to ~\cite{blanchini2008set}. \qed
\end{lemma}
$\Gamma$ can be a design or decision variable. For any positive definite (PD) step matrix $\Gamma,$ there exists a PD value matrix $P$ for which the states of the trajectory-oriented GD \eqref{eq:gdl2} descend to the origin, or equivalently, $\rho (I-2\Gamma P) < 1$, where $\rho$ represents the spectral radius. This is the case where there is no restriction on the behavior of the gradient system \eqref{eq:gdl}. However, our aim is to design a policy $u_k$ to update the states of the system \eqref{eq:sys} using the trajectory-oriented GD dynamics \eqref{eq:gdl} with some step size matrix $\Gamma$. In this case, we must consider the dynamics of \eqref{eq:sys} into account.

\subsection{Design of a static feedback policy}
As mentioned, the idea is to look at the system states $x_k$ as the parameters of the cost function $V_k$ that descend along the gradient of $V_k$. However, in contrast to the usual GD algorithm usage, which is utilized to update the parameters of the controller or the underlying neural network in the opposite direction of a given cost function, we aim to update or shape the states of the system by designing an appropriate control signal $u_k$ using the GD approach. The following theorem summarizes the design procedure. 
\begin{theorem} \label{theo:th1}
    Consider the system \eqref{eq:sys} and assume that the pair $(A, B)$ is stabilizable. Let the control policy be 
    \begin{align}
        u_k = K x_k, \label{eq:uk}
    \end{align}
    where $K \in \mathbb{R}^{m \times n}$. The states of the system \eqref{eq:sys} descend with step size or direction matrix $\Gamma \succ 0$ \eqref{eq:gdl} along the gradient of \eqref{eq:Vk} if there exists $Y \succ 0$ and $F \in \mathbb{R}^{m \times n }$ such that \eqref{eq:gdlmi} holds and 
    \begin{align}
        A Y + B F = Y - 2\Gamma, \label{eq:gdequation}
    \end{align}
    where $F = K Y$. Besides, the closed-loop system is $\lambda-$contractive and parameterized by the triplet of $(\Gamma, P, K)$. Finally, for any stable control gain $K$, there exists a pair of PD $P$ and PD non-symmetric $\Gamma$ to model the closed-loop system as a trajectory-oriented GD system. \\
\textbf{Proof}. According to Lemma \ref{lemma:gd}, the closed-loop system remains stable provided that 
the system \eqref{eq:sys} with control policy $u_k = K x_k$ behaves like the trajectory-oriented GD system \eqref{eq:gdl}. In other words, if
\begin{align}
    A + B K  = I - 2\Gamma P. \label{eq:abkgp}
\end{align}
Multiplying both sides of \eqref{eq:abkgp} by $Y$ results in \eqref{eq:gdequation}. According to \eqref{eq:gdequation}, any stable closed-loop system can be modeled by a trajectory-oriented GD dynamic if a PD step size matrix exists to satisfy it. Based on $I-(A + B K) =  2\Gamma P$, this amounts to showing that $I-(A + B K)$ is PD. This is true since $A+BK$ is a stable matrix. Let $\mu_i$ and $\nu_i$ be the eigenvalues and eigenvectors of $A+BK$ respectively. Then, we have
\begin{align*}
    I \nu_i &= \nu_i, \\
    (A+BK) \nu_i &= \mu_i \nu_i, \\
    (I-(A+BK))\nu_i &= (1-\mu_i) \nu_i.
\end{align*}
So, $1-\mu_i$'s are the eigenvalues of $I-(A+BK)$. Since the closed-loop system $A+BK$ is stable, the real parts of $\mu_i$'s are between $-1$ and $1$, and as a result, the real parts of the eigenvalues of $I-(A+BK)$ are between $0$ and $2$. The step size $\Gamma$ is not symmetric in general because $I-(A+BK)$ is not symmetric in most cases. 
\qed
\end{theorem}

\begin{remark}Note that not all systems can behave in a steepest gradient descent (SGD) manner. In other words, for $\Gamma=\gamma I,$ there is not a static feedback gain $K$ in general to achieve the steepest descent behavior for the closed-loop states.
The fact that any stable closed-loop system can be parameterized by $I-2\Gamma P$ provides an intuitive way to shape the trajectories of a system. In fact, the elements of matrix $\Gamma$ show the weights of the gradient of the cost in shaping the next state of the system. In the context of the equation \eqref{eq:abkgp}, 
$\Gamma$ controls how aggressively or conservatively the feedback gain $K$ is selected in the opposite direction of the cost determined by the matrix $P$.
\end{remark}
 In the next section, we examine the relationship between the trajectory-oriented control using GD (Theorem \ref{theo:th1}) and the well-known LQR approaches. 

% add two theorems here, one for uk general and one for uk=Kxk. 

% \textcolor{red}{any closed loop system can be described by I-gamma P. add lemma
% show that gamma = (I-(A+BK))P-1
% exists for a given K and P. 
% } 

\section{Connection to LQR} \label{sec:4}
Let us consider the infinite-horizon LQR cost function
\begin{align} 
J=\sum_{k=0}^{\infty}{x_k^T Q x_k + u_k^T R u_k},
\label{eq:cost}
\end{align}
which must be minimized by a controller \eqref{eq:uk}. 
The trajectory-oriented control using GD in Theorem \ref{theo:th1} (for $\lambda=1$) has a connection to the LQR problem which is summarized by Lemma \ref{lemma:gd-lqr}.
\begin{lemma} \label{lemma:gd-lqr}
    Let the optimal LQR gain that optimizes \eqref{eq:cost}  be $\bar{K}$. Then, the closed-loop LQR-based control system is equivalent to the trajectory-oriented GD closed-loop system \eqref{eq:gdl2} with 
    \begin{align}
        \bar{\Gamma} =  \frac{1}{2}(I -(A+B \bar{K} )) \bar{P}^{-1}. \label{eq:GammaQR}
    \end{align}
    where $\bar{P}$ is the solution to the LQR Riccati equation \eqref{eq:lqrP}.

\textbf{Proof}. The cost-to-go function for the linear control policy \eqref{eq:uk} can be written as 
\begin{align}
V_k = \min_{K, P} \sum_{t=k}^{\infty}{x_t^T Q x_t + u_t^T R u_t}, \label{eq:VK}
\end{align}
which indeed is equivalent to 
\begin{align}
    V_k &= \min_{K, P} x_k^\top Q x_k + u_k^\top R u_k + V(x_{k+1}). \label{eq:VKVK1}
\end{align}
It is well-known that the value function of LQR can be written as 
\begin{align}
    V_k =  x_k^\top \bar{P} x_k, \quad \bar{P} \succ 0, \label{eq:lqrV}
\end{align}
where
\begin{align}
     \bar{P}=Q+\bar{K}^\top R \bar{K} + (A+B\bar{K})^\top  \bar{P} (A+B \bar{K}), \label{eq:lqrP}
\end{align}
and $\bar{K}$ are the optimal value and grain matrices of the LQR problem, respectively. The LQR-based closed-loop control becomes $x_{k+1}=(A+B \bar K) x_k$. Using \eqref{eq:abkgp}, this closed-loop system corresponds to a trajectory-oriented GD dynamics with an optimal $\bar{\Gamma}$ matrix. That is, 
\begin{align} \label{A+BK}
A+B \bar{K} = I - 2 \bar{\Gamma} \bar{P},
\end{align}
which is used to find \eqref{eq:GammaQR}. Since $(A+B \bar K)$ is stable, the matrix $I-(A+B\bar{K})$ is a non-symmetric PD matrix in general, as shown in Theorem \ref{theo:th1}. Therefore, since $\bar P$ is also PD, for every control gain $\bar K$, a non-symmetric PD matrix $\bar{\Gamma}$ exists to satisfy $I-(A+B\bar{K})=2  \bar{\Gamma} \bar{P}$, and thus \eqref{A+BK}. 
\end{lemma}
The $\bar{\Gamma}$ weights the effect of the LQR value matrix $\bar{P}$ on the states of the optimal closed-loop system, which indeed reflects the effect of $Q$ and $R$ matrices on the trajectories of the system.
Another observation is that the LQR is not a SGD method in the state space because $\bar{\Gamma} \neq \bar{\gamma}I$ in general. This is valid since the SGD is not an optimal way to update the parameters of a learning algorithm, and in our case, the trajectories of the closed-loop system since it does not consider the curvature of the cost function into account. However, letting $\bar{\Gamma}$ to be a PD matrix helps the LQR method to capture more information about the underlying cost function and system dynamics. 
\begin{remark}We can state that for a given PD $\Gamma$, if there exists a solution to Theorem \ref{theo:th1}, it is an optimal $K$ for some $Q$ and $R$. 
\end{remark}
An advantage of the presented approach in Theorem \ref{theo:th1} is that it can be extended to other GD methods to find different forms of trajectory-oriented controllers using the GD concept.
\section{Extension to Heavy-Ball Gradient method} \label{sec:5}
The heavy-ball gradient method (also known as Polyak's heavy-ball method) is an optimization algorithm designed to accelerate GD by adding a momentum term. This method was introduced by Boris Polyak in 1964~\cite{POLYAK19641,lessard2016analysis} and aims to improve the convergence rate of gradient-based optimization algorithms, especially in problems where the objective function is smooth and convex.
The heavy-ball method modifies the traditional GD's update rule by incorporating a momentum term that retains a fraction of the previous step's update. The update rule for the heavy-ball method in terms of the states of the system \eqref{eq:sys} and the cost function \eqref{eq:Vk} can be expressed as
\begin{align}
    x_{k+1} = x_k - 2 \Gamma P x_k + \Delta (x_k - x_{k-1}), \label{eq:hbg}
\end{align}
where the matrix $\Delta$ is the momentum parameter that controls the influence of the previous update direction. The term $\Delta (x_k - x_{k-1})$ adds momentum by incorporating a fraction of the previous step's update to the current step. This momentum term helps the method move faster in directions with consistent gradients, reducing oscillations in directions with varying gradients and can accelerate convergence, especially in optimization landscapes that are ill-conditioned or have long, narrow valleys. The equation \eqref{eq:hbg} can be written in the augmented form $X_k^T = \begin{bmatrix}
    x_k^\top & x_{k-1}^\top
\end{bmatrix}$ as
\begin{align}
    X_{k+1} = \begin{bmatrix}
        I-2\Gamma P + \Delta & -\Delta \\
        I & 0
    \end{bmatrix} X_k. \label{eq:hbgMatrix}
\end{align}
\begin{theorem} \label{lemma:hbgd}
Consider the system \eqref{eq:sys} and assume that the pair $(A,B)$ is stabilizable. Let the control policy be 
    \begin{align}
        u_k = K_1 x_k + K_2 x_{k-1}, \label{eq:ukhb}
    \end{align}
    where $K_1, K_2 \in \mathbb{R}^{m \times n}$. The closed-loop system behaves like \eqref{eq:hbg} with step size matrix or direction matrix $\Gamma \succ 0$ and momentum matrix $\Delta$ if there exists a PD $Y$ and $\Gamma,$  such that 
\begin{subequations}
    \begin{align}
    &\begin{bmatrix}
        \lambda Y & 0 & (Y-2\Gamma+\Delta Y )^\top & Y \\ 
        0 & \lambda Y & -(\Delta Y)^\top & 0 \\
        Y-2\Gamma+\Delta Y   & -\Delta Y & Y & 0 \\
        Y & 0 & 0 &  Y
    \end{bmatrix}  \succeq 0, \label{eq:hbgdlmi} \\
    &A Y + B F_1 = Y-2 \Gamma + \Delta Y, \label{eq:hbk1}\\
    &B F_2 = - \Delta Y, \label{eq:hbk2}
\end{align}
\end{subequations}
hold, where $Y=P^{-1}, F_1 = K_1 Y, F_2=K_2 Y$. Furthermore, the closed-loop system  
\begin{align}
    x_{k+1} = (A+BK_1)x_k + BK_2 x_{k-1}, \label{eq:hbclosed}
\end{align}
is $\lambda-$contractive with $0< \lambda \leq 1$.\\
\textbf{Proof}. Consider the augmented Lyapunov function
\begin{align}
    \bar{V}_k = V_k + V_{k-1} = X_k^\top \begin{bmatrix}
        P & 0\\ 0 & P
    \end{bmatrix} X_k. \label{eq:Vbar}
\end{align}
To hold the $\lambda-$contractivity condition $\bar{V}_{k+1} \leq \lambda \bar{V}_k$ for all $x_k $ and $x_{k-1}$ with a given $0 < \lambda \leq 1$, we can derive the following inequality using the matrix forms of \eqref{eq:hbgMatrix} and \eqref{eq:Vbar}
\begin{dmath}
    \begin{bmatrix}
        I - 2 \Gamma P + \Delta & -\Delta \\ I & 0
    \end{bmatrix}^\top \begin{bmatrix}
        P & 0\\ 0 & P
    \end{bmatrix} \begin{bmatrix}
        I - 2 \Gamma P + \Delta & -\Delta \\ I & 0
    \end{bmatrix}                       \preceq \lambda \begin{bmatrix}
        P & 0\\ 0 & P
    \end{bmatrix}. \label{eq:VKbarineq}
\end{dmath}
By multiplying both sides of \eqref{eq:VKbarineq} with $diag(P^{-1}, P^{-1})$, and using the Schur complement, the inequality \eqref{eq:hbgdlmi} will be extracted. Also, if we compare \eqref{eq:hbclosed} and \eqref{eq:hbg}, the following two equality constraints can be found
\begin{align}
    A  + B K_1 &= I-2 \Gamma P + \Delta, \\
    B K_2 &= - \Delta,
\end{align}
which are respectively equivalents of \eqref{eq:hbk1} and \eqref{eq:hbk2} after multiplying both equality constraints from right by $P^{-1}.$ \qed
\end{theorem}

\section{Simulations} \label{sec:6}
% \textcolor{red}{
% what to be discussed 
% \begin{itemize}
%     \item state trajectories with a scalar step size gamma to show the trajectories that are orthogonal to the level sets
%     \item check lmis with a non symmetric P
%     \item find gamma for a lqr design and compare the results
%     \item compare Nesterov and heavyball results.
% \end{itemize}
% }

\begin{figure}[!t]
\centerline{\includegraphics[width=0.9\columnwidth]{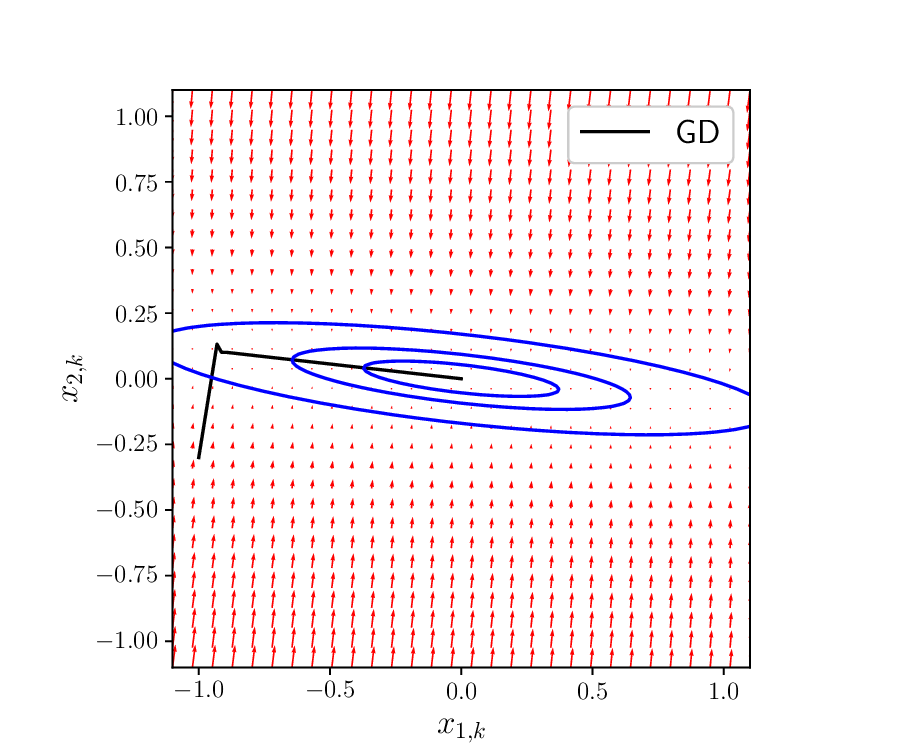}}
\caption{SGD-like behavior of the vehicle steering model for parameters defined in \eqref{eq:gamma}. }
\label{fig:gamma}
\vspace{-14pt}
\end{figure}

\begin{figure}[!t]
\centerline{\includegraphics[width=0.9\columnwidth]{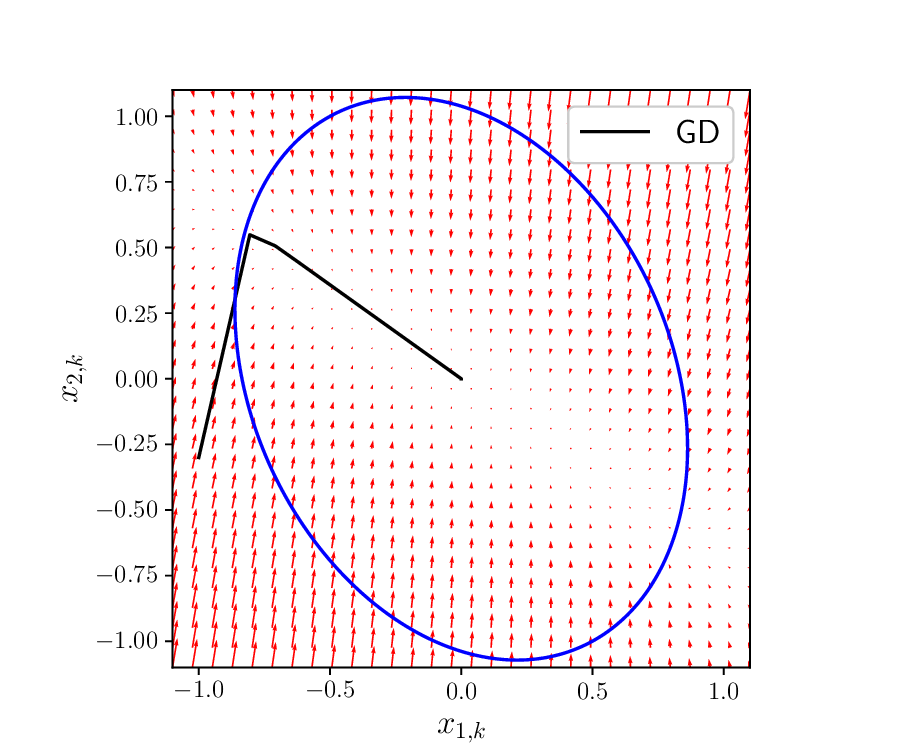}}
\caption{GD-like behavior of the vehicle steering model for parameters defined in \eqref{eq:Gamma} }
\label{fig:Gamma}
\vspace{-15pt}
\end{figure}

Consider the vehicle steering model (or lateral deviation dynamics)~\cite{kishida2022risk} with the LTI model \eqref{eq:sys} where the lateral position and heading angle are $x_{1,k}$ and $x_{2,k}$, respectively. The system matrices are
$
    A=\begin{bmatrix}
        1 & 0.2 \\ 0 & 1
    \end{bmatrix}, B=\begin{bmatrix}
        0.06 \\ 0.2
    \end{bmatrix}.
$
First, we want to show the behavior of the controlled system when the $\Gamma = \gamma I.$ Using Theorem \ref{theo:th1}, we obtain
\begin{align}
    P = \begin{bmatrix}
        1.69 & 5.65 \\ 5.65 & 32.95
    \end{bmatrix}, K = \begin{bmatrix}
        -1.33 & -7.76
    \end{bmatrix}, \gamma = 0.0015. \label{eq:gamma}
\end{align}
% \gamma = 0.02357
Figure \ref{fig:gamma} shows the trajectory of the vehicle steering system for the initial condition $x_0=\begin{bmatrix}
    -1& -0.3
\end{bmatrix}^\top.$ The system's states evolve in the SGD manner. It is also evident that the trajectory flows in the opposite direction of the cost function gradient as the trajectory is perpendicular to the level-sets of $V_k$ (blue ellipsoids.) It is unclear how we can achieve this trajectory using an LQR controller. 
For a non-symmetric PD $\Gamma$ (Figure \ref{fig:Gamma}), the solution to Theorem \ref{theo:th1} is
\begin{align}
    P = \begin{bmatrix}
        1.43 & 0.29 \\ 0.29 & 0.93
    \end{bmatrix}, K = \begin{bmatrix}
        2.86 & -4.59
    \end{bmatrix}, \Gamma = \begin{bmatrix}
        0.055 & 0.023 \\ 0.107 & 0.462
    \end{bmatrix}. \label{eq:Gamma}
\end{align}
% Figure \ref{fig:Gamma} shows the trajectory of the system using parameters in \eqref{eq:Gamma}. 

Now, let us compare the performance of the trajectory-oriented controller in Theorem \ref{theo:th1} to the LQR with $Q=I$ and $R=10.$ The LQR value and gain matrices are respectively
\begin{align}
    P = \begin{bmatrix}
        13.08 & 13.30 \\ 13.30 & 37.63
    \end{bmatrix}\,\, \text{and} \,\, K = \begin{bmatrix}
        -0.29 & -0.76
    \end{bmatrix}. \label{eq:lqrdesign}
\end{align}
According to \eqref{eq:GammaQR}, $\Gamma=\begin{bmatrix}
    0.0085 & -0.0071 \\ 0.00052 & 0.0038
\end{bmatrix}$  corresponds to the LQR design \eqref{eq:lqrdesign}. The trajectory-oriented controller using GD for this $\Gamma$ (Theorem \ref{theo:th1}) is
\begin{align}
    P = \begin{bmatrix}
        5.13 & 5.38 \\ 5.38 & 17.59
    \end{bmatrix}, K = \begin{bmatrix}
        -0.23 & -0.69
    \end{bmatrix}, \label{eq:gd_lqr}
\end{align}
% Figure \ref{fig:gd-lqr} shows the states and controls obtained using the LQR and GD methods. 
% Gamma: 
%  [[ 0.0005     -0.00170518]
%  [ 0.00142893  0.00198728]]
The small difference between the gain in \eqref{eq:gd_lqr} and \eqref{eq:gd_lqr} is due to the numerical solution for the LMI.
The trajectory-oriented control using GD method results in an optimal method in the sense of the LQR for the given values of $Q$ and $R.$ The cost associated with the LQR \eqref{eq:lqrdesign} and GD \eqref{eq:gd_lqr} designs are $24.45$ and $24.99$ respectively. 

\begin{figure}[!t]
\centerline{\includegraphics[width=0.9\columnwidth]{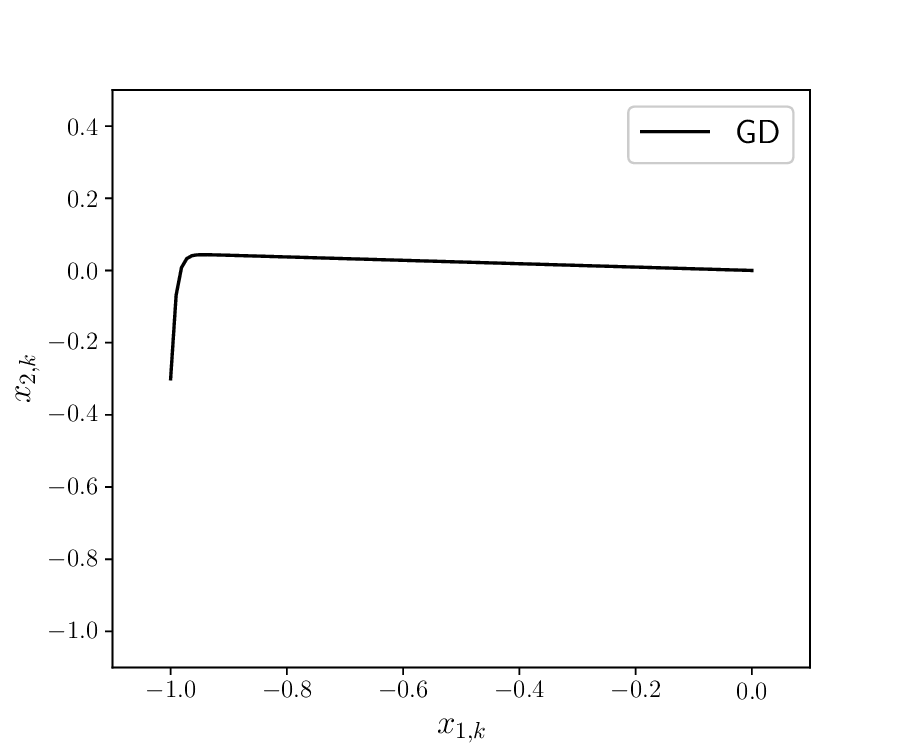}}
\caption{Explicitly shaping the behavior of the vehicle steering model using $\Gamma$}
\label{fig:x2}
\vspace{-10pt}
\end{figure}

\indent A significant advantage of our approach is that we can leverage the direction matrix $\Gamma$ to shape the trajectories of the closed-loop system. Let us add the constraint $\Gamma \preceq \begin{bmatrix}
        0.001 & 0 \\ 0 & 0.1
    \end{bmatrix}$ to Theorem \ref{theo:th1}.
This condition limits the convergence speed of the first state compared to the second one. As a result, if the dynamics of the system allows, the second state converges to the origin faster than the first, as shown in Figure \ref{fig:x2}. Here are the design results
\begin{align}
     P &= \begin{bmatrix}
        6.18 & 0.19 \\ 0.19 & 4.72
    \end{bmatrix}, K = \begin{bmatrix}
        -0.15 & -3.34
    \end{bmatrix}, \nonumber \\
    \Gamma &= \begin{bmatrix}
        0.00075 & 0.00004 \\ 0.00032 & 0.071
    \end{bmatrix}. \label{eq:x2-dom}
\end{align}
% \begin{figure}[!t]
% % \centerline{\includegraphics[width=\columnwidth]{images/lqr-gd-traj.eps}}
% \centerline{\includegraphics[width=\columnwidth]{images/lqr-gd-x1.eps}} 
% \centerline{\includegraphics[width=\columnwidth]{images/lqr-gd-x2.eps}}
% \centerline{\includegraphics[width=\columnwidth]{images/lqr-gd-uk.eps}}
% \caption{Gradient descent behavior of the vehicle steering model for parameters defined in \eqref{eq:Gamma} }
% \label{fig:lqr-gd}
% \end{figure}

% \begin{figure}[htbp] % 'htbp' specifies figure placement (here, top, bottom, or on a separate page)
%   \centering
%   % Subfigure (a)
%   \begin{subfigure}[b]{0.45\textwidth} % Width of subfigure (adjust as needed)
%       \centering
%       \includegraphics[width=\textwidth]{images/lqr-gd-x1.eps} % Replace with your image file
%       % \caption{$x_{1,k}$}
%       \label{fig:subfig1}
%   \end{subfigure}
%   \hfill % Space between subfigures

%   % Subfigure (b)
%   \begin{subfigure}[b]{0.45\textwidth}
%       \centering
%       \includegraphics[width=\textwidth]{images/lqr-gd-x2.eps}
%       % \caption{$x_{2,k}$}
%       \label{fig:subfig2}
%   \end{subfigure}

%   % Subfigure (c)
%   \begin{subfigure}[b]{0.45\textwidth}
%       \centering
%       \includegraphics[width=\textwidth]{images/lqr-gd-uk.eps}
%       % \caption{$u_k$}
%       \label{fig:subfig3}
%   \end{subfigure}
  
%   \caption{Comparing states and controls of the LQR \eqref{eq:lqrdesign} and GD \eqref{eq:gd_lqr} designs}
%   \label{fig:gd-lqr}
% \end{figure}

\section{Conclusion} \label{sec:7}
The proposed method offers a different perspective on controller design by employing the GD algorithm to directly update the states of a dynamical system. Unlike traditional approaches that focus on optimizing control parameters or minimizing cost functions indirectly, our method shapes the system's state trajectories by descending along the gradient of a state-dependent cost function. This approach introduces a new degree of flexibility and adaptability in control design, allowing for more precise tuning of system behavior via the step size matrix $\Gamma$. Furthermore, by establishing a connection between the trajectory-oriented GD-based controller and the well-known LQR, we show that our method not only broadens the scope of gradient-based optimization in control applications but also provides a novel framework for understanding and designing optimal controllers in both linear and nonlinear settings.  The future work will investigate the similarities and differences between the presented trajectory-oriented control using the GD method and classical control approaches such as model reference control, pole placement, and eigenstructure assignment. Their similarities and differences will then be leveraged to shape the trajectories generated by any of these controllers to correct their undesired behaviors. 

\bibliographystyle{IEEEtran}
\bibliography{refs} 
\end{document}